\documentclass[reqno, 11pt]{amsart}

\usepackage{color}    
\usepackage{mathtools} 
\usepackage{amssymb}

\makeatletter
\def\ckech{\mathaccent"\accentclass@014}
\def\tah{\mathaccent"\accentclass@05E}
\makeatother

  \renewcommand\check{\bm\ckech} 
\renewcommand\hat{\bm\tah}

\newcommand\mysection[1]{%
\section{#1}\setcounter{equation}{0}}

\makeindex

\usepackage{bm}

\newtheorem{theorem}{Theorem}[section]

\newtheorem{corollary}[theorem]{Corollary}

\theoremstyle{definition}

\theoremstyle{remark}

\newtheorem{remark}[theorem]{Remark}

\newcommand{\vsharp}{\asymp\kern -.5em\|}

 \makeatletter 
 \def\dashint{\operatorname{\,\,\,\mathclap{\!\int}\! \!\text{\bf--}\!\!}}  
 \makeatother

\newcommand\bC{\mathbb{C}}

\newcommand\bR{\mathbb{R}}
\newcommand\bQ{\mathbb{Q}}
\newcommand\bS{\mathbb{S}}

\newcommand\cL{\mathcal{L}}

\newcommand\cO{\mathcal{O}}
\newcommand\cP{\mathcal{P}}

  {{}}

\renewcommand\){{\rm)}}

\def\+){\tmspace+\thinmuskip{.05em}\)}

\def\dashnorm{\,\,\text{\bf--}\kern-.5em\|}

\renewcommand{\eqref}[1]{\text{(\ref{#1})}}

\begin{document}

\title[]{Extending BMO functions
in parabolic setting}
\author[]{N.V. Krylov}
\address{School of Mathematics, University of Minnesota, Minneapolis, MN, 55455}
\email{nkrylov@umn.edu}



\begin{abstract}
We prove that one can extend any $BMO^{x}$
function $a$ given in a cube in $\bR^{d+1}$
to become a $BMO^{x}$ functions $\hat a$ in $\bR^{d+1}$
almost preserving its $[a]^{\sharp}$ seminorm,
which is, loosely speaking, $L_{\infty}$-norm
of the maximal function 
in $t$ and $BMO$-norm in $x$.

\end{abstract}
\maketitle
\mysection{Introduction}

This article is devoted to proving
that one can extend functions,
given in a parabolic box and
being measurable in $t$ and BMO
in $x$, to the whole space in
a ``controlled'' way, so that
the $BMO^{x}$ characteristic
of the extended function
is controlled by the $BMO^{x}$ characteristic
of the original one.
Our interest arose from the necessity
to extend matrix-valued functions
(second-order coefficient in parabolic
operators). In addition, we need to
preserve the nondegeneracy and boundedness
of extended functions.
 
In the case of functions independent of $t$
the problem of extending BMO functions
given in a domain
to the whole space, with controlled
increase of the BMO seminorm, received
much attention. The first paper settling
the issue belongs to P. Jones, see  \cite{Jo_80},
where he puts forward very powerful theorem
with necessary and sufficient conditions
on the domain for such extension to be possible. Later (see
\cite{CW_15} and the references therin)
it turned out that John's conditions on the domain mean exactly that it should be a uniform domain (see the definition,
for instance, in \cite{CW_15}, where the authors mention that balls are uniform
domains). P. Jones dealt with bounded
domains, the case of unbounded ones
is treated, for instance, in \cite{BD_23}.

Recall that if $\cO$ is a domain in $\bR^{d}=\{x=(x^{1},...,x^{d}):x^{i}\in\bR\}$ and $u$ is a function on $\cO$,
then we write $u\in BMO(\cO)$ if for every
cube $Q=(a,b)^{d}\subset \cO$ the quantity
\begin{equation}
                            \label{7.10.1}
\frac{1}{|Q|^{2}}\int_{Q}\int_{Q}
|u(x)-u(y)|\,dxdy\quad (|Q|=\text{vol}\,Q)
\end{equation}
is bounded by a constant independent of $Q$.
The least such constant is denoted by
$[u]_{BMO(\cO)}$. If the quantity
\eqref{7.10.1} goes to zero as $|Q|
\downarrow 0$ we write $u\in VMO(\cO)$.

Up to now in most applications of BMO functions
in PDEs one does not treat so general
$\cO$ as in \cite{Jo_80}, \cite{CW_15},
\cite{BD_23}, and for simple domains checking the conditions in these papers requires 
diving into a deep area of mathematics
which is quite far from PDEs. Furthermore,  sometimes
we need extensions to have additional properties such as vanishing at certain
distance from the original domain.
Therefore, for instance, in \cite{Ac_92}
the author proves from scratch that
VMO functions in bounded
 domains with Lipschitz boundary can be extended  to 
VMO functions in $\bR^{d}$ and 
in such a way that the extension
vanishes outside   a slightly larger domain.
It is worth saying that his approach does not work for extending
BMO functions.

In \cite{LOR_21} the following theorem
is proved without using sophisticated
terminology from \cite{Jo_80}, \cite{CW_15},
\cite{BD_23} (the statement is slightly changed).
\begin{theorem}
                     \label{theorem 7.10.2}
Let $A_{1}=(-1,1)^{d}, A_{2}=(-2,2)^{d}$,
and let $a\in BMO(A_{2})$. Then there is
a function $\hat a$ on $\bR^{d}$ such that
\begin{equation}
                       \label{7.10.4}
\hat a=a\quad\text{on}\quad A_{1},\quad
\hat a=\frac{1}{|A_{2}|}\int_{A_{2}}a\,dx\quad\text{outside}
\quad A_{2},
\end{equation}
\begin{equation}
                       \label{7.10.5}
[\hat a]_{BMO(\bR^{d})}\leq N[a]_{BMO(A_{2})},
\end{equation}
where $N$ is independent of $a$.
\end{theorem}

\begin{corollary}
                  \label{corollary 7.10.1}
For $\delta\in(0,1]$ define $\bS_{\delta}$
as the set of $d\times d$-symmetric matrices
whose eigenvalues are in $[\delta,\delta^{-1}]$ and assume that we are given
an $\bS_{\delta}$-valued function
$a\in BMO(A_{2})$. Then there is an $\bS_{\delta}$-valued function
$\hat a$ on $\bR^{d}$ such that \eqref{7.10.4}
and \eqref{7.10.5} hold.

\end{corollary}

Indeed, in the first place, this follows from the proof
in \cite{LOR_21}. But also let us call $\check a^{ij}$ the function
we obtain by applying Theorem \ref{theorem 7.10.2} to $a^{ij}$ and then define
$\hat a$ as the projection of $\check a=
(\check a^{ij},i,j=1,...,d)$ as a point in
the space of all $d\times d$-matrices
on $\bS_{\delta}$. Then note that if $\Pi$
is the projection operator,  then
$|\Pi(\check a(x))-\Pi(\check a(y))|
\leq |\check a(x)-\check a(y)|$.

\mysection{Main result}

Unlike Corollary \ref{corollary 7.10.1}
here we are dealing with $\bS_{\delta}$-valued functions $a(t,x)$, $(t,x)\in
\bR^{d+1}=\bR\times\bR^{d}$, which are measurable in $t$ and kind of BMO in  $x$.
For $x\in\bR^{d}$ set
$$
 \|x\|=\max_{i\leq d}|x^{i}|.
$$
If
$\Gamma\subset\bR^{d+1} $ we denote
by $|\Gamma|$ its Lebesgue measure and when
it makes sense for functions $f(t,x)$ we
denote
$$
\dashint_{\Gamma}f\,dxdt=\frac{1}{|\Gamma|}
\int_{\Gamma}f\,dxdt.
$$

 For $\rho>0$ the collection of
$$
x+(0,\rho)^{d},\quad x\in\bR^{d},
$$
is denoted by $\bC_{\rho}$. The collection of
$$
(s,s+\rho^{2})\times C,\quad s\in\bR, C
\in\bC_{\rho},
$$
is denoted
by $\bQ_{\rho}$. If $Q=(a,b)\times C\in 
\bQ_{\rho}$ we set
$$
\rho_{Q}=\rho,\quad Q^{z}=(a,b),\quad Q^{\tau}=C.
$$
Set
$$
\cP=(0,1)^{d+1}.
$$
 
  For domains $\cO\subset\bR^{d+1}$ we denote by $\bQ_{\rho}(\cO)$
the collection of $Q\in \bQ_{\rho}$
such that $Q\subset \cO$.

For a domain
$\cO\subset\bR^{d+1}$,
$\rho_{0}\in(0,\infty]$, and a function $a$ given in $\cO$
  we write $a\in BMO^{x}_{\rho_{0}}(\cO)$ if
$$
[a]^{\sharp}_{\rho_{0},\cO }:=\sup_{\rho\leq \rho_{0}}
\sup_{ Q\in\bQ_{\rho} (\cO)}
\dashint_{ Q ^{z}}
\dashint_{Q ^{\tau}}\dashint_{Q ^{\tau}}|a (t,x)-a (t,y)|\,dxdydt<\infty.
$$

\begin{remark}
                     \label{remark 7.10.1}
Observe that if $a$ is indpendent of $t$
and $a\in BMO^{x}_{1}(\cP)$, then
$a\in BMO((0,1)^{d})$. However, if $a$ is independent of $x$, then $[a]^{\sharp}_{\rho_{0},\cO }=0$. This somewhat justifies
the notation $BMO^{x}$.
\end{remark}

\begin{theorem}
                    \label{theorem 7.8.1}
Let $a$ be a  real-valued
function of class
$ BMO^{x}_{1}(\cP)$.
Then there exists an $\bS_{\delta}$-valued
function $\hat a$ on $\bR^{d+1}$ 
such that  
\begin{equation}
                         \label{7.11.1}
a=\hat a\quad\text{on}\quad\cP,\quad
[\hat a]^{\sharp}_{\infty,\bR^{d+1} }
\leq 3^{2+2d} [a]^{\sharp}_{1,\cP } .
\end{equation}
 
\end{theorem}

The proof of the theorem is based
on the following particular case of the change of variables formula,
Theorem 2 in Section 3.3.3 of \cite{EG_15}.

\begin{theorem}
                     \label{theorem 7.11.2}
Let $\phi:\bR^{d}\to\bR^{d}$ have finite
Lipschitz constant and $\hat a\in L_{1}(\bR^{d})$. Denote  by $J$ the Jacobian of
$\phi$. Then
\begin{equation}
                        \label{7.11.30}
\int_{\bR^{d}}\hat a J\,dx=
\int_{\bR^{d}}\Big(\sum_{x:\phi(x)=y}\hat a(x)\Big)\,dy
\end{equation}

\end{theorem}
\begin{remark}
                      \label{remark 7.11.2}
If $\hat a(x)=a(\phi(x))$, equation \eqref{7.11.30} becomes
\begin{equation}
                        \label{7.11.3}
\int_{\bR^{d}}a(\phi(x)) J(x)\,dx=
\int_{\bR^{d}}a(y)\#\{x:\phi(x)=y\} \,dy.
\end{equation}
This formula is proved if $\hat a\in L_{1}(\bR^{d})$. For obvious reasons it is also true for any $a\geq0$.
\end{remark}

{\bf Proof of Theorem \ref{theorem 7.8.1}}. Set $a(t,x)=1$
if $(t,x)\in \partial \cP$.  Introduce $2$-periodic function
$\phi(t)$, $t\in\bR$, such that $\phi(t)=t$
on $[0,1]$ and $\phi(t)=1-t$ on $[1,2]$, and set
$$
\phi( x)=\big( \phi(x^{1}),...,\phi(x^{d})\big),\quad
\hat a(t,x)=a\big(\phi(t),\phi( x))\big).
$$
Observe that the Lipschitz 
constant of the mapping
$$
(t,x,y)\to \big(\phi(t),\phi(x),\phi(y)\big) 
$$ is equal to one
and its  Jacobian equals one
(a.e.). Therefore, by Remark~\ref{remark 7.11.2}, for any $Q\in\bQ$
$$
I_{Q}:=\int_{Q^{z}}\int_{Q^{\tau}}\int_{Q^{\tau}}
|\hat a(t,x)-\hat a(t,y)|\,dxdydt
$$
$$
=\int_{\bR^{2d+1}}|a(\phi(t),\phi(x))
-a(\phi(t),\phi(y))|I_{Q^{z}}(t)I_{Q^{\tau}}
(x)I_{Q^{\tau}}(y)\,dxdydt
$$
\begin{equation}
                            \label{7.9.2}
=\int_{(0,1)^{2d+1}}|a(t_{1},x_{1})-a(t_{1},y_{1})|\sigma^{z}_{Q}(t_{1})\sigma^{\tau}_{Q}(x_{1})\sigma^{\tau}_{Q}(y_{1})\,dx_{1}dy_{1}dt_{1},
\end{equation}
 where $\sigma_{Q}(t_{1})$
is the number of points in
$$
 \{t\in Q^{z}
 :
  \phi(t) = t_{1}  \}
$$
and $\sigma_{Q}(x_{1})$ 
is the number of points in
$$
 \{x\in Q^{\tau}
 :
  \phi(x) = x_{1}  \}.
$$

{\em Case $\rho_{Q}\leq1$\/}. In that
case 
\begin{equation}
                          \label{7.11.4}
\sigma_{Q}(t_{1})\leq 2, \quad \sigma_{Q}(x_{1}),\sigma_{Q}(y_{1})\leq 2^{d}.
\end{equation}
Indeed, given time interval $[c,c+1)$,
each value of $\phi$ on this interval can
repeat itself at most twice. The same
is true for the number of $x^{i}\in [c,c+1)$
for which $\phi(x^{i})$ takes a fixed
value. Since there are $d$ of the $x^{i}$'s,
$\sigma_{Q}(x_{1})\leq 2^{d}$.

Furthermore, if $\sigma^{z}_{Q}(t'_{1})\ne0$ and $\sigma^{z}_{Q}(t''_{1})\ne0$, then
for some $t',t''\in Q^{z}$ we have
$$
|t_{1}'-t_{1}''|=|\phi(t')-\phi(t'')|\leq |t'-t''|\leq\rho_{Q}^{2}.
$$
It follows that there exists $t_{0}\in
(0,1)$ such that
$$
\{t_{1}\in(0,1):\sigma^{z}_{Q}(t_{1})\ne0\}
\subset [t_{0}-\rho_{Q}^{2}/2,t_{0}+\rho_{Q}^{2}/2]\subset[0,1].
$$
Similarly, if $\sigma^{\tau}_{Q}(x'_{1})\ne0$ and $\sigma^{\tau}_{Q}(x''_{1})\ne0$, then
for some $x',x''\in Q^{\tau}$ we have
$$
\|x_{1}'-x_{1}''\|=\|\phi(x')-\phi(x'')\| \leq\rho_{Q} .
$$
We conclude that there exists $\hat C\in \bC_{\rho_{Q}}$, $\hat C
\subset [0,1]^{d+1}$ such that
$$
\{x_{1}\in(0,1)^{d}:\sigma^{\tau}_{Q}(x_{1})\ne0\}
\subset \hat C,
$$
  For
$$
\hat Q=(t_{0}-\rho_{Q}^{2}/2,t_{0}+\rho_{Q}^{2}/2\)\times \hat C
$$
this yields $\hat Q\in \bQ_{\rho_{Q}}(\cP)$,
$$
I_{Q}\leq 2^{2d+1}\int_{\hat Q^{z}}\int_{(0,1)^{d}}\int_{(0,1)^{d}}|a(t ,x )-a(t ,y )| I_{\hat Q^{\tau}}(x)I_{\hat Q^{\tau}}(y)\,dx dy dt  
$$
$$
= 2^{2d+1}\int_{\hat Q^{z}}\int_{\hat Q^{\tau}}\int_{\hat Q^{\tau}}|a(t ,x )-a(t ,y )|  \,dx dy dt\leq 2^{2d+1}\rho^{2d+2}_{Q}a^{\sharp}_{1,\cP}.
$$

{\em Case $\rho_{Q}\in(n-1,n]$ with an integer $n\geq 2$\/}. In this case in formula \eqref{7.9.2} we have
$$
\sigma^{z}\leq n^{2}+1<(n+1)^{2},\quad\sigma^{\tau}
\leq (n+1)^{d}
$$
and then
$$
I\leq  (n+1)^{2+2d}
\int_{(0,1)}\int_{(0,1)^{d}}\int_{(0,1)^{d}}|a(t ,x )-a(t ,y )|  \,dx dy dt.
$$
The integral term above is less than 
$a^{\sharp}_{1,\cP}$ by definition
and 
$$
(n+1)^{2+2d}\leq \rho_{Q}^{2+2d}\big(\frac{n+1}{n-1}\Big)^{2+2d}\leq \rho_{Q}^{2+2d}
3 ^{2+2d}.
$$
Therefore,
$$
I\leq 3^{2+2d} a^{\sharp}_{1}\rho^{2+2d}_{Q}. 
$$
This proves the theorem. \qed

By applying the above proof
to matrix-valued functions, we get
the following.

\begin{corollary}
                   \label{corollary 7.11.1}

Let $a$ be an $\bS_{\delta}$-valued
function of class
$ BMO^{x}_{1}(\cP)$.
Then there exists an $\bS_{\delta}$-valued
function $\hat a$ on $\bR^{d+1}$ 
such that \eqref{7.11.1} holds.
\end{corollary}

\begin{remark}
                   \label{remark 7.10.2}
From the above proof one can see that if
$[a]^{\sharp}_{\rho,\cP }\to 0$ as
$\rho\downarrow 0$, then also
$[\hat a]^{\sharp}_{\rho,\bR^{d+1} }\to 0$  as
$\rho\downarrow 0$.

\end{remark}
\begin{remark}
                   \label{remark 7.11.1}
If $a$ is independent of $t$,
we obtain an extension theorem
for usual BMO functions.

\end{remark}
\begin{remark}
                   \label{remark 7.11.3}
A typical situation in which 
Corollary \ref{corollary 7.11.1}
is useful is the following.
Let $a$ be an $\bS_{\delta}$-valued
function of class
$ BMO^{x}_{\infty}(\bR^{d+1})$,
$$
\cL u=\partial_{t}u+a^{ij}D_{ij}u.
$$
It is relatively easy (see, for instance,
\cite{Kr_26}) to prove that, if $p\in(1,\infty)$ and
\begin{equation}
                         \label{7.11.6}
[a]^{\sharp}_{\infty,\bR^{d+1}}\leq \check a(d,\delta,p)
\end{equation}
for an appropriate $\check a(d,\delta,p)>0$, then for any $\lambda\geq0$ and $u\in C^{\infty}_{0}(\bR^{d+1})$
\begin{equation}
                         \label{7.11.5}
\|\partial_{t}u,D^{2}u,\sqrt\lambda Du,
\lambda u\|_{L_{p}}\leq N(d,\delta,p)
\|\cL u-\lambda u\|_{L_{p}},
\end{equation}
where $L_{p}=L_{p}(\bR^{d+1})$.
Below by $N$ we denote generic constants
depending only on $d,\delta,p$.

The global condition \eqref{7.11.6} is
inconvenient in many instances especially
if the ultimate goal is to consider equations in bounded domains.

Let us replace it with
\begin{equation}
                         \label{7.11.7}
[a]^{\sharp}_{1,(t,x)+\cP}\leq \check a(d,\delta,p)\quad \forall (t,x)\in\bR^{d+1}.
\end{equation}

It that case for $u\in C^{\infty}_{0}(\bR^{d+1})$ which has support in $\cP$ we have
$$
\hat \cL u:=\partial_{t}u+\hat a^{ij}D_{ij}u
=\cL u
$$
on $\bR^{d+1}$ and, consequently,
\eqref{7.11.5} holds. This is true
if $u$ has support in any of $(t,x)+\cP$.
Then take a nonnegative $\zeta\in C^{\infty}_{0}(\cP)$ such that the integral
of its $p$th power is one and introduce
$$
\zeta_{s,y}(t,x)=\zeta (t+s,x+y),\quad u_{s,y}=u\zeta_{s,y}.
$$
By the above, for instance,
$$
\int_{\bR^{d+1}}|D^{2}u_{s,y}|^{p}\,dxdt
\leq N\int_{\bR^{d+1}}|\cL u_{s,y}-\lambda
u_{s,y}|^{p}\,dxdt.
$$

Here
$$
D_{ij}u_{s,y}=\zeta_{s,y}D_{ij}u+
D_{i}\zeta_{s,y}D_{j}u+D_{j}\zeta_{s,y}D_{i}u+uD_{ij}\zeta_{s,y},
$$
$$
2^{-p}\zeta^{p}_{s,y}
|D^{2}u|^{p}\leq |D^{2}u_{s,y}|^{p}  +N|D\zeta_{s,y}|^{p}\,|Du|^{p}
 +N|D^{2}\zeta_{s,y}|^{p}|u|^{p}
$$
(we used that $(a+b)^{p}\leq 2^{p}(a^{p}+b^{p}), a,b\geq0$). Similarly.
$$
2^{-p}\zeta^{p}_{s,y}
|\partial_{t}u|^{p}\leq |\partial_{t}u_{s,y}|^{p}  +N|\partial_{t}\zeta_{s,y}|^{p}\,| u|^{p} ,
$$
$$
2^{-p}\zeta^{p}_{s,y}
|Du|^{p}\leq |Du_{s,y}|^{p}  +N|D\zeta_{s,y}|^{p}\,| u|^{p} .
$$

Futhermore,
$$
|\cL u_{s,y}-\zeta_{s,y}\cL u|\leq
N\big(|u|\,|\cL \zeta_{s,y} |+|D\zeta_{s,y}| \,|Du|\big).
$$
Hence,
$$
\|\zeta_{s,y}\partial_{t}u,\zeta_{s,y}D^{2}u,\zeta_{s,y}\sqrt\lambda Du,
\zeta_{s,y}\lambda u\|^{p}_{L_{p}}\leq N(d,\delta,p)
\|\zeta_{s,y}(\cL u-\lambda u)\|^{p}_{L_{p}}
$$
$$
+N\|uD^{2}\zeta_{s,y}|^{p}_{L_{p}}+N\|u\partial_{t}\zeta_{s,y}|^{p}_{L_{p}}
+\| N|D\zeta_{s,y}|Du\|
^{p}_{L_{p}}.
$$
After integrating through this relation
with respect to $(s,y)\in\bR^{d+1}$
and observing that, for instance,
$$
\int_{\bR^{d+1}}\| \zeta_{s,y}D^{2}u \|^{p}_{L_{p}}\,dyds
$$
$$
=\int_{\bR^{d+1}}|D^{2}u(t,x)|^{p}
\Big(\int_{\bR^{d+1}}\zeta_{s,y}^{p}(t,x)
\,dyds\Big)\,dxdt=\| D^{2}u \|^{p}_{L_{p}},
$$
we come to the following
$$
\|\partial_{t}u,D^{2}u,\sqrt\lambda Du,
\lambda u\|_{L_{p}}\leq N 
\|\cL u-\lambda u\|_{L_{p}}+N_{1}\|Du\|_{L_{p}}+N_{0}\|u\|_{L_{p}},
$$
which is valid for all $\lambda\geq0$.
One can absorb the terms on the right with $Du$ and $u$
into the left-hand side if
$$
2N_{0}\leq\lambda,\quad 2N_{1}\leq \sqrt\lambda.
$$
In that case one obtains \eqref{7.11.5}
under condition \eqref{7.11.7} for all
$\lambda\geq \lambda_{0}=\lambda_{0}(d,\delta,p)>0$ and any $u\in C^{\infty}_{0}
(\bR^{d+1})$.
\end{remark}

{\bf Acknowledgement}. The author is sincerely grateful to A. Lerner for attracting his attention to \cite{LOR_21}
and fruitful discussion.


\begin{thebibliography}{m}

\bibitem{Ac_92}
P. Acquistapace, {\em On BMO regularity for linear elliptic systems\/}, Ann. Mat. Pura Appl. (4), Vol. 161 (1992), 231--269.

\bibitem{BD_23} A. Butaev and G. Dafni,  {\em Locally uniform domains and extension of bmo functions\/}, Annales Fennici Mathematici, Vol. 48 (2023),
No. 2, 567--594.

\bibitem{CW_15}
Yuming Chu and Gendi Wang, {\em Two necessary and sufficient conditions for uniform domains\/}, Hokkaido Math. J., Vol. 35 (2006), No. 4, 935--942. 

\bibitem{EG_15}
L.C. Evans and R.F. Gariepy,  Measure theory and fine properties of functions. Revised edition. Textbooks in Mathematics. CRC Press, Boca Raton, FL, 2015. xiv+299 pp.

\bibitem{Jo_80}
P.W. Jones, {\em Extension theorems for BMO\/}, Indiana Univ. Math. J., Vol. 29 (1980),
No. 1, 41--66.

\bibitem{Kr_26} N.V. Krylov,
{\em Essentials of Real Analysis and Morrey-Sobolev spaces for second-order elliptic and parabolic PDEs with singular first-order coefficients\/}, arXiv:2505.14863

\bibitem{LOR_21} 
A.K. Lerner, S. Ombrosi, I.P. Rivera-R\'ios, {\em On two weight estimates for iterated commutators\/}, J. Funct. Anal., Vol. 281 (2021), No. 8, Paper No. 109153, 46 pp.

\end{thebibliography}
\end{document}